\documentclass{article}

\usepackage{amsmath}
\usepackage{amssymb}
\usepackage{amsthm}
\usepackage[all,cmtip,knot]{xy}
\usepackage{graphicx}

\theoremstyle{plain} \newtheorem{proposition}[]{Proposition}[section]
\theoremstyle{plain} \newtheorem{lemma}[]{Lemma}[section]
\theoremstyle{plain} \newtheorem{theorem}[]{Theorem}[section]
\theoremstyle{plain} \newtheorem{corollary}[]{Corollary}[section]
\theoremstyle{plain} \newtheorem{conjecture}[]{Conjecture}[section]
\theoremstyle{definition} \newtheorem{definition}[]{Definition}[section]
\theoremstyle{remark} \newtheorem{remark}[]{Remark}[section]

\newcommand{\proj}{\mathcal{P}}
\newcommand{\destab}{\mathcal{D}}
\newcommand\Z{\mathbb{Z}}
\newcommand\Q{\mathbb{Q}}
\def\fin\qedhere
\def\pr {{\text{pr}}}
\DeclareMathOperator{\Hom}{Hom}

\date{}

\title{On the Infinity Flavor of Heegaard Floer Homology and the Integral Cohomology Ring}
\author{Tye Lidman}

\begin{document}
\maketitle

\begin{abstract}
In \cite{hfpa}, Ozsv\'ath and Szab\'o construct a spectral sequence with $E_2$ term $\Lambda^*(H^1(Y;\mathbb{Z}))\otimes \mathbb{Z}[U,U^{-1}]$ converging to $HF^\infty(Y,\mathfrak{s})$ for a torsion Spin$^c$ structure $\mathfrak{s}$.  They conjecture that the differentials are completely determined by the integral triple cup product form via a proposed formula given in \cite{plumbed}.  In this paper, we prove that $HF^\infty(Y,\mathfrak{s})$ is in fact determined by the integral cohomology ring when $\mathfrak{s}$ is torsion.  Furthermore, for torsion Spin$^c$ structures, we give a complete calculation of $HF^\infty$ with mod 2 coefficients when $b_1$ is 3 or 4.
\end{abstract}

\section{Introduction}

Throughout the previous decade, Heegaard Floer theory has been a very useful and calculable machine in low-dimensional topology.  It includes invariants for closed three- and four-manifolds, as well as for knots and links.  Similarly, manifolds with boundary, singular knots, contact structures, and other objects can be studied as well.  One of the most effective computational tools in Heegaard Floer theory is the integral surgery formula \cite{hfkz}, which converts the Heegaard Floer complexes of a closed, oriented 3-manifold $Y$ and a nullhomologous knot $K \subset Y$ into the Heegaard Floer homology of surgeries on $K$.
Given a Heegaard splitting of $Y$ along a surface $\Sigma$, Heegaard Floer homology calculates the Lagrangian Floer homology of tori in the symmetric product of $\Sigma$.  The Heegaard Floer homology of $Y$ splits as a direct sum over the set of Spin$^c$ structures on $Y$.  Different flavors of Heegaard Floer homology twist the differential by a count of the intersection number of a holomorphic disk with some choice of basepoint(s) on the surface.

While many new results in low-dimensional topology have come from calculations of these groups, one flavor, $HF^\infty$, has the simplest structure.  In fact, it has been calculated for $b_1(Y) \leq 2$ in \cite{hfpa}.  In this case, it is completely determined by the integral cohomology ring.  Also, Mark \cite{thommark} has obtained results in this direction, gaining information about $HF^\infty$ from a complex $C_*^\infty(Y)$ with differential given completely by the triple cup product.  If one calculates $HF^\infty$ with coefficients in  $\mathbb{Z}[[U,U^{-1}]$ instead, it is shown in \cite{hflz} that these groups, $\bold{HF}^\infty(Y,\mathfrak{s})$, vanish for any non-torsion Spin$^c$ structure $\mathfrak{s}$.  Therefore, we are only concerned with torsion Spin$^c$ structures in this paper.  From now on, any Spin$^c$ structure $\mathfrak{s}$ will be torsion.
In \cite{hfpa}, it is shown that for each torsion Spin$^c$ structure $\mathfrak{s}$ there exists a spectral sequence with $E_2$ term $\Lambda^*(H^1(Y;\mathbb{Z})) \otimes \mathbb{Z}[U,U^{-1}]$ converging to $HF^\infty(Y,\mathfrak{s})$.  Furthermore, in \cite{plumbed}, Ozsv\'ath and Szab\'o propose:

\begin{conjecture}
\label{theconjecture}
The differential $d_3 : \Lambda^i(H^1(Y;\mathbb{Z})) \otimes U^j \rightarrow \Lambda^{i-3}(H^1(Y;\mathbb{Z})) \otimes U^{j-1}$ is given by
\begin{equation} \label{differential}
\phi^1 \wedge \ldots \wedge \phi^i \mapsto \frac{1}{3!(i-3)!}\sum_{\sigma \in S_i} (-1)^{|\sigma|} \langle\phi^{\sigma(1)} \smile \phi^{\sigma(2)} \smile \phi^{\sigma(3)},[Y]\rangle \cdot \phi^{\sigma(4)} \wedge \ldots \wedge \phi^{\sigma(i)}
\end{equation}
Furthermore, all higher differentials vanish.  (For notational purposes, we will omit the $U$'s in the domain and range from now on).
\end{conjecture}

Note that if this conjecture is true, knowing the integral triple cup product form on $Y$ allows a complete calculation of $HF^\infty(Y,\mathfrak{s})$.
The goal of this paper is to present a few partial results in this direction.

\begin{theorem}
\label{cohomologydependence}
For $\mathfrak{s}$ torsion, $HF^\infty(Y,\mathfrak{s})$ is completely determined by the integral cohomology ring.  In other words, if $H^*(Y_1;\mathbb{Z}) \cong H^*(Y_2;\mathbb{Z})$ as graded rings and $\mathfrak{s}_1$ and $\mathfrak{s}_2$ are torsion \emph{Spin}$^c$ structures on $Y_1$ and $Y_2$ respectively, then $HF^\infty(Y_1,\mathfrak{s}_1;\mathbb{Z}) \cong HF^\infty(Y_2,\mathfrak{s}_2;\mathbb{Z})$.
\end{theorem}

Often, it will be useful to use coefficients for $HF^\infty$ in $\mathbb{F} = \mathbb{Z}/2$ as opposed to $\mathbb{Z}$.  For notational purposes, when referring to Conjecture~\ref{theconjecture}, we will be taking the integral triple cup product and then reducing mod 2 in the case of $\mathbb{F}$ coefficients.

\begin{theorem}
\label{b1=3}
Conjecture~\ref{theconjecture} holds for $b_1(Y)=3$ with coefficients in $\mathbb{F}$.
\end{theorem}

\begin{theorem}
\label{b1=4}
For $b_1(Y)=4$, $HF^\infty(Y,\mathfrak{s};\mathbb{F})$ agrees with the prediction for the homology given by the conjecture.
\end{theorem}

We now outline the arguments given for the proofs in this paper.  In order to calculate $HF^\infty(Y)$ in general, we prove that it suffices to consider any manifold which can be obtained from $Y$ by a sequence of nonzero surgeries on nullhomologous knots.  This is done by showing that such a surgery does not affect the integral triple cup product form or $HF^\infty$.  Furthermore, we show that we only need to calculate $HF^\infty$ in the case of $H_1(Y;\mathbb{Z}) \cong \mathbb{Z}^n$, by showing that in each torsion Spin$^c$ structure, $HF^\infty(Y,\mathfrak{s})$ behaves as $HF^\infty$ of a manifold which is some ``torsionless'' version of $Y$.

We then use a theorem of Cochran, Geiges, and Orr \cite{surgeryequivalence} which generalizes Casson's result that any integral homology sphere can be related by a sequence of $\pm 1$-surgeries on nullhomologous knots in $S^3$.  This shows that there exists a nice class of ``model manifolds''.  This collection has the property that given any $Y$, there exists some model manifold which can be related to $Y$ by a sequence of such surgeries, and thus has isomorphic $HF^\infty$.  For $b_1 = 3$ and 4, we can explicitly write down these models and calculate $HF^\infty$ simply based on knowledge of $HF^\infty(\mathbb{T}^3,\mathfrak{s}_0)$ and the integral surgery formula of \cite{hfkz}. \\ \\
\textbf{Acknowledgements:} I would like to thank Ciprian Manolescu for his knowledge and patience as an advisor, as well as for sharing with me his construction of homologically split surgery presentations.  I would also like to thank Liam Watson for his encouragement to work on this problem and his aid in drawing Heegaard diagrams.

\section{Eliminating Torsion} \label{eliminatingtorsion}
The goal of this section is to reduce the calculation to the case where $H_1(Y;\mathbb{Z})$ is torsion-free.  The idea is to construct a sufficiently nice surgery presentation and then argue that we can remove each knot that is not contributing to $b_1(Y)$ without changing either the integral triple cup product form or $HF^\infty$.

Let us first start with an example.  Choose a link $L$ in $S^3$ and let $Y$ be the result of surgery on $L$ with framing $\Lambda$.  Fix an integer $n \neq 0$ and perform $n$-surgery on $K$, a knot in $S^3$ separated from $L$ by an embedded sphere.  The resulting manifold will be $Y \# S^3_n(K)$.  Notice that the integral triple cup product form of $Y \# S^3_n(K)$ is isomorphic to that of $Y$.  Similarly, the connect-sum formula for $HF^\infty$ and the calculation of $HF^\infty$ for rational homology spheres of \cite{hfpa} give $HF^\infty(Y \# S^3_n(K), \mathfrak{s}_Y \# \mathfrak{s}_K) \cong HF^\infty(Y,\mathfrak{s}_Y)$ for any choice of Spin$^c$ structures on $Y$ and $S^3_n(K)$.  We have now, in a sense, removed the $n$-surgery on $K$ from $Y \# S^3_n(K)$, and thus removed a factor of $\mathbb{Z}/n$ from $H_1$, but preserved the integral triple product form and $HF^\infty$.

We want to generalize and repeat this procedure in order to remove all of the knots contributing to torsion in $b_1$.

\begin{proposition}
\label{triplecupsurgerypreserved}
Perform $n$-surgery on a nullhomologous knot $K$ in $Y$ for some nonzero integer $n$.  The resulting manifold, $Y_n(K)$, has the same integral triple cup product form as $Y$.
\end{proposition}
\begin{proof}
Using results of Cochran, Gerges, and Orr on rational surgery equivalence \cite{surgeryequivalence}, it suffices to show that there is a sequence of non-longitudinal surgeries on rationally nullhomologous knots, beginning with $Y$ which terminate at $Y_n(K)$.  This is clear since $n \neq 0$.
\end{proof}

The following proposition is noted in Section 4.1 of \cite{plumbed}.

\begin{proposition}[Ozsv\'ath-Szab\'o]
\label{hfsurgerypreserved}
Let $s_K$ be a torsion Spin$^c$ structure on $Y_n(K)$ agreeing with $s$ on $Y - K$.  The resulting manifold has $HF^\infty(Y_n(K),\mathfrak{s}_K) \cong HF^\infty(Y, \mathfrak{s})$.
\end{proposition}

To remove torsion from $H_1$, we need a sufficiently nice surgery presentation to try to repeat the argument of the example.  However, since a surgery presentation might not consist of all pairwise-split components, we have to find the closest thing.  The idea is to represent $Y$ by surgery on a link where the components have pairwise linking number 0.  Such a link is called \emph{homologically split}.  The following theorem tells us that we can do this if we are willing to slightly change the manifold.  The proof can be found at the end of this paper.

\begin{lemma}[Manolescu]
\label{homsplit}
Let $Y$ be a closed, oriented 3-manifold.  There exist finitely many nonzero integers, $m_1, \ldots, m_k$, such that there exists a homologically split surgery presentation for $Y \# L(m_1,1) \# \ldots \# L(m_k,1)$.
\end{lemma}

Now, begin with a homologically split link presentation for $Y'$, the result of $Y$ after summing with the necessary lens spaces.  Since each of the nonzero surgeries will now be performed on a nullhomologous knot, Proposition~\ref{triplecupsurgerypreserved} and
Proposition~\ref{hfsurgerypreserved} show it suffices to calculate $HF^\infty$ of the 3-manifold obtained by surgery on the sublink consisting of components that are 0-framed.  However, if $L$ is a homologically split link with $l$ components, then after 0-surgery on each component, the resulting three-manifold will have $H_1(Y) \cong \mathbb{Z}^l$.
As connect sums with lens spaces do not change $HF^\infty$ or the integral triple cup product, we are content to work in the case that $H_1(Y;\mathbb{Z})$ is torsion-free.

While knowing $HF^\infty$ might not tell us everything we would want to know about the spectral sequence, we will later see that for $b_1 = 3$, this does give complete information.

\section{Model Manifolds}

Following \cite{surgeryequivalence}, we will call two 3-manifolds, $Y_1$ and $Y_2$, \emph{surgery equivalent} if there is a finite sequence of $\pm1$ surgeries on nullhomologous knots, beginning with $Y_1$ and terminating at $Y_2$.

\begin{theorem}[Cochran-Gerges-Orr]
\label{torsionfreecupequivalent}
Let $H_1(Y_1;\mathbb{Z}) \cong \mathbb{Z}^n$.  Suppose that $Y_1$ and $Y_2$ have isomorphic integral triple cup product forms.  Then $Y_1$ and $Y_2$ are surgery equivalent.
\end{theorem}

It is important to note that this is not true if $H_1$ has torsion.  A counterexample can be exhibited by taking $Y_1$ as $\#_{i=1}^3 L(5,1)$ and $Y_2$ as 5-surgery on each component of the Borromean rings (Example 3.15 of \cite{surgeryequivalence}).  However, since each of these has $b_1 = 0$, we know they must have isomorphic $HF^\infty$.  Therefore, this invariant cannot quite detect the subtlety seen by the cup products with other coefficient rings.

\begin{proof}[Proof of Theorem~\ref{cohomologydependence}]
Theorem~\ref{torsionfreecupequivalent} combined with the work of Section~\ref{eliminatingtorsion} now proves that the integral triple cup product determines $HF^\infty$ (for any coefficients).  A little more work allows the statement for the integral cohomology ring.
If the integral cohomology rings of $Y_1$ and $Y_2$, are isomorphic (grading preserving), then the integral triple cup product form of $Y_1$ is isomorphic to either that of $Y_2$ or $-Y_2$.  Note that if we follow the method of removing torsion from $Y_1$ and $Y_2$, then the resulting manifolds, $Y_1'$ and $Y_2'$, will have isomorphic cohomology rings.  Thus, since this action does not change $HF^\infty$ for a given Spin$^c$ structure, we may assume the manifolds do not have torsion in $H_1$.  If $Y_1$ and $Y_2$ have isomorphic forms, then we are clearly done by the theorem.  On the other hand, if $-Y_2$, then we apply \cite{surgeryequivalence} (Corollary 3.8), to see that $Y_2$ is surgery equivalent to $-Y_2$.  This completes the proof.
\end{proof}

As there is an explicit way to construct a 3-manifold with $H_1(Y) = \mathbb{Z}^n$ and arbitrary triple cup product form by essentially iterating ``Borromean surgeries" (see Corollary 3.5 in \cite{surgeryequivalence}), $HF^\infty$ should be very accessible.  In fact, we speculate that the methods of this paper can be generalized to the link surgery formula of Manolescu and Ozsv\'ath \cite{hflz} to calculate the homology for all $b_1$ via this approach.

In the cases $b_1=3$ and $b_1=4$, we can explicitly see what the set of surgery equivalence classes is that we are dealing with.  The following is calculated in Example 3.3 in \cite{surgeryequivalence}.

\begin{theorem}[Cochran-Gerges-Orr]
The surgery equivalence classes for three-manifolds with $H_1(Y;\mathbb{Z}) =
\mathbb{Z}^3$ are precisely determined by $|H^3(Y)/(H^1(Y) \smile
H^1(Y) \smile H^1(Y))|$.  A representative of each is given by
the manifold $M_n$, with Kirby diagram given by \\

\qquad\quad\includegraphics[scale=.5]{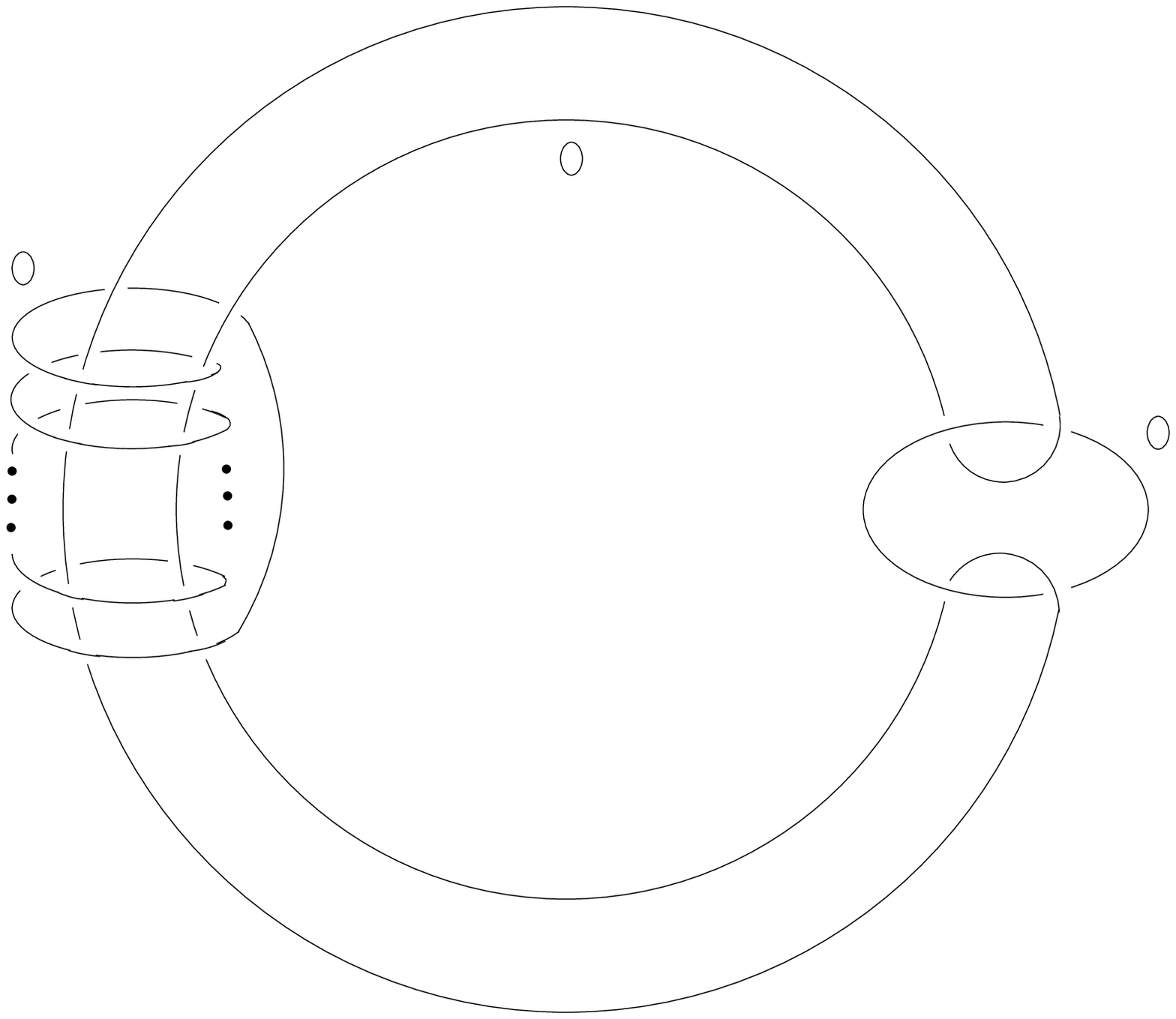}
\begin{center}
\emph{Figure 3.1.}
\end{center}

\end{theorem}

\qquad

It is useful to note that $M_0 = \#_{i=1}^3 S^2 \times S^1$ and $M_1
= \mathbb{T}^3$.  Calculating $HF^\infty$ for this class of
manifolds is what suffices to prove Theorem~\ref{b1=3}. Furthermore,
it turns out that calculating $b_1=3$ combined with the connect-sum
formula is sufficient to understand $b_1 = 4$ as well.

\begin{proposition}[Cochran-Gerges-Orr]
\label{foursplits}
If $H_1(Y) = \mathbb{Z}^4$, then $Y$ is surgery equivalent to $M_n \# S^2 \times S^1$ for some $n \geq 0$.
\end{proposition}

\section{Review of the Surgery Formula}
In this section we review the knot surgery formula from \cite{hfkz} with the perspective and notation of \cite{hflz}.  For the rest of the paper, we will assume that our manifold $Y$ has $H_1(Y;\mathbb{Z})$ torsion-free and that all Heegaard Floer coefficients are $\mathbb{F}$.  Furthermore, we will assume all diagrams are admissible as needed and we are working over the unique torsion Spin$^c$ structure, $\mathfrak{s}_0$, on the relevant manifold.  Let $K$ be a nullhomologous knot in $Y$.  Knowledge of the knot Floer complex will be used to calculate the Heegaard Floer homology of surgeries on $K$.  First, we must introduce the necessary definitions.

Let $(\Sigma,\alpha,\beta,z,w)$ be a Heegaard diagram for $K$ in $Y$.  Note that $(\Sigma,\alpha,\beta,z)$ and $(\Sigma,\alpha,\beta,w)$ are each diagrams for $Y$, and thus no longer contain any information about the knot.  There exists an Alexander grading on $\mathbb{T}_\alpha \cap \mathbb{T}_\beta$ satisfying
\begin{equation}
A(x) - A(y) = n_z(\phi) - n_w(\phi)
\end{equation}
for $\phi \in \pi_2(x,y)$, which can be canonically made into an absolute grading.  Similarly, for any pointed Heegaard diagram for $Y$, there is an absolute Maslov grading (since $\mathfrak{s}_0$ is torsion) satisfying
\begin{equation}
M(x) - M(y) = \mu(\phi) - 2n_p(y),
\end{equation}
where $p$ is the chosen basepoint and again $\phi \in \pi_2(x,y)$.  Recall that multiplication by $U$ lowers $A$ by 1 and $M$ by 2.
We can now define a $CFK$-like complex with differential twisted by the Alexander grading.  Let $x \vee y = \max\{x,y\}$.

\begin{definition}
 $\mathfrak{A}_0$ is the chain complex over $\mathbb{F}[U,U^{-1}]$ freely-generated by the subset of $\mathbb{T}_\alpha \cap \mathbb{T}_\beta$ consisting of elements corresponding to $\mathfrak{s}_0$ equipped with the differential
\begin{equation}
\partial_0 x = \sum_{y \in \mathbb{T}_\alpha \cap \mathbb{T}_\beta} \sum_{ \phi \in \pi_2(x,y), \mu(\phi) = 1 } \#(\mathcal{M}(\phi)/\mathbb{R}) \cdot U^{A(x) \vee 0 - A(y) \vee 0 + n_w(\phi)} y
\end{equation}
for $x \in \mathbb{T}_\alpha \cap \mathbb{T}_\beta$.
\end{definition}

There exist chain maps relating $\mathfrak{A}_0$, $CF_z = CF^\infty(\Sigma,\alpha,\beta,z,\mathfrak{s}_0)$, and $CF_w = CF^\infty(\Sigma,\alpha,\beta,w,\mathfrak{s}_0)$ as given by

\[
\xymatrix{
& \mathfrak{A}_0 \ar[ld]_{\proj^{-K}} \ar [rd]^{\proj^K}\\
CF_z \ar[rr]^{\destab^{-K}}& & CF_w}
\]
where the diagonal maps are the projections $\proj^K(x) = U^{A(x) \vee 0}x$ and $\proj^{-K}(x) = U^{(-A(x)) \vee 0}x$.  After possible stabilizations of the diagram avoiding both $z$ and $w$, the diagrams $(\Sigma,\alpha,\beta,z)$ and $(\Sigma,\alpha,\beta,w)$ can be related by a sequence of isotopies and handleslides, since they both represent $Y$.  Choose such a sequence of moves and let the corresponding induced composition of triangle-counting maps on $CF^\infty$ be denoted $\destab^{-K}$.  Note that this map is unique up to chain homotopy.  We can assume $\destab^K$ to be the identity map, coming from a choice of no isotopies or handleslides from $(\Sigma,\alpha,\beta,w)$ to itself.

\begin{proposition} \label{gradingpreserved}
 The projection maps, $\proj^{\pm K}$ are quasi-isomorphisms which preserve relative Maslov grading.  Furthermore, $\destab^{-K}$ is a quasi-isomorphism and its induced map on homology preserves absolute gradings.
\end{proposition}
\begin{proof}
The projections are quasi-isomorphisms simply because they are defined by $U$ multiplications and are therefore bijective chain maps.  That they preserve the relative grading follows from the work on gradings in \cite{hflz}.  On the other hand, $\destab^{-K}$ consists of a composition of triangle counts, each of which is a chain homotopy equivalence that preserves the absolute grading on $HF$ \cite{absgraded}.
\end{proof}

Following \cite{hflz}, let $\Phi^{-K} = \destab^{-K} \circ \proj^{-K}$ and $\Phi^K = \destab^K \circ \proj^K = \proj^K$.  We now state the surgery formula of Oszv\'ath and Szab\'o for $0$-surgery on $K$.

\begin{theorem}[Ozsv\'ath-Szab\'o] \label{surgeryformula}
Suppose $\mathfrak{s}_0$ is torsion.  Consider the chain map $\Psi^K = \Phi^{-K} + \Phi^K$ from $\mathfrak{A}_0$ to $CF_w$.  The homology of the mapping cone of $\Psi^K$ is isomorphic to $HF^\infty(Y_0(K),\mathfrak{s})$, where $\mathfrak{s}$ is the unique torsion Spin$^c$ structure that agrees with $\mathfrak{s}_0$ on $Y-K$.
\end{theorem}

For notational convenience, let $H_*(CF_p) = \mathcal{K}_p$ for $p = z$ or $w$, and $H_*(\mathfrak{A}_0) = \mathcal{K}_{z,w}$.  It is important to note that from the surgery formula, $HF^\infty(Y_0(K),\mathfrak{s})$ is an $\mathbb{F}[U,U^{-1}]$-module of rank equal to rk$\mathcal{K}_w+$ rk$\mathcal{K}_{z,w}-$ $2$rk$(\Psi^K_*)$.
We will abuse notation and use the same symbols for chain maps and their induced maps on homology when there is minimal confusion.

\section{Example: $\mathbb{T}^3$} \label{examplesection}
Recall that we are interested in calculating the Heegaard Floer homology of the manifolds $M_n$ given in Figure 3.1.  The main goal of this section is to understand the simplest nontrivial example, $M_1 = \mathbb{T}^3$.  From Figure 3.1, we can represent $M_n$ by 0-surgery on a knot in $S^2 \times S^1 \# S^2 \times S^1$ and can therefore apply the surgery formula.  For $M_1$, this in fact gives 0-surgery on the Borromean rings.  Let $K$ denote the remaining knot in 0-surgery on two components of the Borromean rings.  Although the Heegaard Floer homology of $\mathbb{T}^3$ has already been calculated in \cite{absgraded}, reverse-engineering the computation via the surgery formula will allow us to deduce valuable information for the general case.  First, let us study the differentials of the spectral sequence for $b_1 = 3$.

From (\ref{differential}), it is clear that after $d_3$ the spectral sequence must collapse.  In fact, the only possibly nontrivial component of $d_3$ maps from $\Lambda^3(H^1)$ to $\Lambda^0(H^1)$, each of which has rank 1.  Therefore, to prove the result for $b_1 = 3$, it suffices to find $HF^\infty$.  If it has rank 8, then $d_3 \equiv 0$, and if it has rank 6, then $d_3(\phi^1 \wedge \phi^2 \wedge \phi^3) = 1$.
Before dealing with $M_1$, we note that $M_0 = \#_{i=1}^3 S^2 \times S^1$ has $HF^\infty(M_0) \cong \mathbb{F}[U,U^{-1}] \otimes H^*(\mathbb{T}^3;\mathbb{F})$, which corresponds to $d_3$ being identically 0 in (\ref{differential}).
For $\mathbb{T}^3$, Conjecture~\ref{theconjecture} predicts that the map $d_3:\Lambda^3 (H^1) \rightarrow \Lambda^0(H^1)$ should be nonzero, which agrees with rk$HF^\infty(\mathbb{T}^3,\mathfrak{s}_0) = 6$.  We now want to use this to understand the map $\destab^{-K}$ in detail.  The best way to see this is via matrix representations, so we must pick out the right bases for $\mathcal{K}_{z,w}, \mathcal{K}_z,$ and $\mathcal{K}_w$.
Since $\mathcal{K}_{z,w} \cong \mathcal{K}_z \cong \mathcal{K}_w \cong \mathbb{F}[U,U^{-1}] \otimes H^*(\mathbb{T}^2)$, we can choose two $\mathbb{F}$-bases $x_1, x_2$ and $y_1,y_2$ for $\mathcal{K}_z$ at adjacent Maslov gradings (say $x_1,x_2$ for $(\mathcal{K}_z)_0$ and $y_1,y_2$ for $(\mathcal{K}_z)_1$).  This clearly gives an ordered  $\mathbb{F}[U,U^{-1}]$-basis for the entire module.

Define the map $\Theta^K:CF_z \rightarrow CF_w$ by $\Theta^K(x)=U^{A(x)}x$.  Like the projection maps, this is a quasi-isomorphism.

\begin{proposition} \label{projectioncomposetheta}
$\Theta^K \circ \proj^{-K} = \proj^K$.
\end{proposition}
\begin{proof}
Add the powers of $U$ together.
\end{proof}

\begin{lemma}  $\Theta^K$ preserves relative Maslov grading.  Furthermore, it preserves the parity of the absolute Maslov grading.  (We will see later that it preserves absolute grading after more work).
\end{lemma}
\begin{proof} The first statement simply follows from Proposition~\ref{projectioncomposetheta} and Proposition~\ref{gradingpreserved}.
For the second statement, we must rely on the spectral sequence.  In order for the rank of $HF^\infty$ to be at least 6, the rank of the matrix representation for $\Phi^K + \Phi^{-K}$ on homology must be 0 or 1.  Since $\Phi^K + \Phi^{-K} = (\Theta^K + \destab^{-K}) \circ \proj^{-K}$, factoring out the the quasi-isomorphism $\proj^{-K}$ shows $\Theta^K + \destab^{-K}$ must also have rank 0 or 1.  Suppose that $\Theta^K$ reverses the parity of the absolute grading.  Choose an ordered basis for $\mathcal{K}_w$ as two pairs of elements in adjacent Maslov gradings as before, where the parities agree with the ordering for $\mathcal{K}_z$.  Since $\destab^{-K}$ preserves the grading parity and $\Theta^K$ reverses it, the maps are represented by:
\begin{equation*}
\destab^{-K} =
\begin{pmatrix}
A & 0 \\
0 & B
\end{pmatrix} \quad \text{and} \quad
\Theta^K =
\begin{pmatrix}
0 & C \\
D & 0
\end{pmatrix},
\end{equation*}
where $A,B,C,D \in GL_2(\mathbb{F}[U,U^{-1}])$.  This implies $\destab^{-K} + \Theta^K$ has rank at least 2, which is a contradiction.
\end{proof}

Since $\Theta_K$ preserves the relative grading and its parity, we can choose a basis for $\mathcal{K}_w$ such that $\Theta_K$ is represented by $U^k$ times the identity, for some $k \in \mathbb{Z}$.  Necessarily, this basis will be supported in a pair of adjacent Maslov gradings and ordered such that the parities agree with $\mathcal{K}_z$.
Note that $\destab^{-K}$ is now represented by a matrix ($x_1,x_2,y_1,y_2$ is the ordering) of the form

\begin{equation*}
\begin{pmatrix}
a & b & 0 & 0 \\
c & d & 0 & 0 \\
0 & 0 & e & f \\
0 & 0 & g & h
\end{pmatrix} \qquad\quad a,b,c,d,e,f,g,h \in \mathbb{F}.
\end{equation*}

Choose a basis for $\mathcal{K}_{z,w}$ such that $\proj^{-K}$ can be represented by the identity.  Thus, the last thing that we want to understand is the matrix representation of $\proj^K$.

\begin{lemma}
With respect to these bases, $\proj^K$ and $\Theta^K$ are both represented by the identity.
\end{lemma}
\begin{proof}
Because the representation for $\proj^{-K}$ is the identity, Proposition~\ref{projectioncomposetheta} guarantees $\proj^K$ and $\Theta^K$ will be represented by the same matrix.
If $\proj^K = U^k I$ for some $k \neq 0$, then
\begin{equation*}
\Phi^K + \Phi^{-K}=
\begin{pmatrix}
U^k + a & b & 0 & 0 \\
c & U^k + d & 0 & 0 \\
0 & 0 & U^k + e & f\\
0 & 0 & g & U^k + h
\end{pmatrix}
\end{equation*}
must have rank at least 2.  However, from our previous remarks about the rank of the homology of the mapping cone, this would contradict the rank of $HF^\infty(\mathbb{T}^3,\mathfrak{s}_0)$ being 6.
\end{proof}

Consider the collection of matrices

\begin{equation*}
X = \left\{
\begin{pmatrix}
1 & 1 & 0 & 0 \\
0 & 1 & 0 & 0 \\
0 & 0 & 1 & 0 \\
0 & 0 & 0 & 1
\end{pmatrix},
\begin{pmatrix}
1 & 0 & 0 & 0 \\
1 & 1 & 0 & 0 \\
0 & 0 & 1 & 0 \\
0 & 0 & 0 & 1
\end{pmatrix},
\begin{pmatrix}
0 & 1 & 0 & 0 \\
1 & 0 & 0 & 0 \\
0 & 0 & 1 & 0 \\
0 & 0 & 0 & 1
\end{pmatrix}, \right.
\end{equation*}
\begin{equation*}
\qquad \qquad \qquad \left. \begin{pmatrix}
1 & 0 & 0 & 0 \\
0 & 1 & 0 & 0 \\
0 & 0 & 1 & 1 \\
0 & 0 & 0 & 1
\end{pmatrix},
\begin{pmatrix}
1 & 1 & 0 & 0 \\
0 & 1 & 0 & 0 \\
0 & 0 & 1 & 0 \\
0 & 0 & 1 & 1
\end{pmatrix},
\begin{pmatrix}
1 & 1 & 0 & 0 \\
0 & 1 & 0 & 0 \\
0 & 0 & 0 & 1 \\
0 & 0 & 1 & 0
\end{pmatrix} \right\}
\end{equation*}

We now see these matrices in fact describe the only possibilities for $\destab^{-K}$.

\begin{proposition} \label{matrixpossibilities}
The map $\destab^{-K}$ is represented by a matrix in $X$.
\end{proposition}

\begin{proof}
Here, we explicitly use the fact that the rank of $\Phi^K + \Phi^{-K}$ must be precisely 1.  This is because $HF^\infty(\mathbb{T}^3,\mathfrak{s}_0)$ has rank 6 and both $\mathcal{H}_w$ and $\mathcal{H}_{z,w}$ have rank 4.  Since $\Phi^K + \Phi^{-K}$ is represented by
\begin{equation*}
\begin{pmatrix}
a+1 & b & 0 & 0 \\
c & d+1 & 0 & 0 \\
0 & 0 & e+1 & f \\
0 & 0 & g & h+1
\end{pmatrix},
\end{equation*}
exactly three of the two-by-two blocks must be identically 0 and the other must have rank 1.  It is easy to check that each of the matrices in $X$ have this property.  We see that either $\begin{pmatrix} a & b \\ c & d\end{pmatrix}$ or $\begin{pmatrix} e & f \\ g & h\end{pmatrix}$ is the identity.  Without loss of generality, we assume
$\begin{pmatrix}
e & f \\
g & h
\end{pmatrix} =
\begin{pmatrix}
1 & 0 \\
0 & 1
\end{pmatrix}.$
Now, the possible matrices
$\begin{pmatrix}
a & b \\
c & d
\end{pmatrix} \in GL_2(\mathbb{F})$
that don't show up in $X$ are
$\begin{pmatrix}
1 & 0 \\
0 & 1
\end{pmatrix},
\begin{pmatrix}
0 & 1 \\
1 & 1
\end{pmatrix}$, and
$\begin{pmatrix}
1 & 1 \\
1 & 0
\end{pmatrix}.$
Direct calculation shows that $\Phi^K + \Phi^{-K}$ would have either rank 0 or rank 2 in any of these cases, which would be a contradiction.  Repeating the arguments with the top-left and bottom-right blocks switched discounts all of the other matrices not in $X$.
\end{proof}

We have now obtained the necessary info, namely $\destab^{-K}$, to generalize.  The important thing to note is that for 0-surgery on any nullhomologous knot in $Y$, where $H_1(Y;\mathbb{Z}) = \mathbb{Z}^2$, all of the arguments other than that of Proposition~\ref{matrixpossibilities} are still valid, since the only piece of information used was that the spectral sequence guarantees the rank of $HF^\infty$ is at least 6.

However, Proposition~\ref{matrixpossibilities} does not apply to 0-surgery on every knot in $S^2 \times S^1 \#S^2 \times S^1$.  Doing 0-surgery on a split unknot, $K_0$, to get $\#_{i=1}^3 S^2 \times S^1$ , which has rank 8, shows that $\Phi^{K_0} = \Phi^{-K_0}$.  This in fact means that after these choice of bases, $\destab^{-K_0}$ must be the identity.  Therefore, following this framework, it remains to calculate $\destab^{-Z_n}$ to obtain $HF^\infty$ for the other knots, $Z_n$, yielding $M_n$ ($n\geq 2$).
To do this, we iterate the $\mathbb{T}^3$ calculation repeatedly using a technique we call composing knots.

\section{Composing Knots and the Calculation for $M_n$}
Recall that given a Heegaard diagram $(\Sigma,\alpha,\beta)$, any two points on $\Sigma - \alpha - \beta$ determine a knot, $K$, in $Y$.  Now, suppose there are instead 3 distinct points, $z$, $u$, and $w$.  Then the pairs of basepoints, $(z,u), (u,w), (z,w)$, determine three knots.  We want to consider Heegaard diagrams containing this information.  Also, we will ignore orientations as this will not affect $HF^\infty$ or the predicted differentials in the spectral sequence.

\begin{definition}
A \emph{Heegaard diagram} for $(K_1,K_2,K)$ in $Y$ is a Heegaard diagram for $Y$, $(\Sigma, \alpha, \beta)$, equipped with 3 distinct basepoints $z$, $u$, and $w$, in $\Sigma - \alpha - \beta$, such that $(z,u)$,$(u,w)$, and $(z,w)$ determine $K_1$,$K_2$, and $K$ respectively.
\end{definition}

\begin{proposition} \label{destabscompose}
Consider a Heegaard diagram for $(K,K_1,K_2)$.  At the chain level, $\destab^{-K} = \destab^{-K_2} \circ \destab^{-K_1}$ up to chain homotopy.  Hence, the equality holds on homology.
\end{proposition}
\begin{proof}
$\destab^{-K_1}$ is induced by a sequence of Heegaard moves taking $(\Sigma,\alpha,\beta,z)$ to $(\Sigma,\alpha,\beta,u)$ and $\destab^{-K_2}$ comes from a sequence of moves from $(\Sigma,\alpha,\beta,u)$ to $(\Sigma,\alpha,\beta,w)$.  Therefore, the composition of isotopies and handleslides goes from $(\Sigma,\alpha,\beta,z)$ to $(\Sigma,\alpha,\beta,w)$ and induces a map $\destab^{-K}$.
\end{proof}

Thus, it is important to note that since most of the complexity in the knot surgery formula comes from the map $\destab^{-K}$, having a Heegaard diagram for $(K,K_1,K_2)$ and an understanding of each $\destab^{-K_i}$ should make the computation more manageable.  This is the approach we will use for the rest of the $M_n$.  However, we must first establish that such things exist and more importantly, derive a way of relating this information to the $M_n$.

\begin{lemma} \label{morsehelmet}
Suppose $K_1$ and $K_2$ are knots in $Y$ where $K_1 \cap K_2$ is an embedded connected interval.  Then if $K = (K_1 \cup K_2) - K_1 \cap K_2$, there exists a Heegaard diagram for $(K,K_1,K_2)$.  (See below)
\end{lemma}

\xygraph
{
  !{0;/r14.0pc/:}
  !{\hcap-|{\displaystyle{K_1-K_2}}}
  !{0;/r7.0pc/:}
  [r]!{\xcapv[2]@(0)|{\displaystyle{K_1 \cap K_2}}}
  !{0;/r14.0pc/:}
  !{\hcap|{\displaystyle{K_2-K_1}}}
}
\begin{center}
Figure 6.1. Each simple cycle corresponds to a knot.
\end{center}

\begin{proof}
The idea follows the construction of Heegaard diagrams for knots in \cite{introhf}.  Begin with a self-indexing Morse function, $h:S^3 \rightarrow [0,3]$, with exactly two critical points.  Note that traversing a flow from index 0 to index 3 and then another in ``reverse'' gives a knot.  Thus, three flow lines give three knots in a natural way as before.

\[
\xymatrix{
3 \ar@/_2pc/[dd] \ar[dd] \ar@/^2pc/[dd]\\ \\
0
}
\]
\begin{center}
Figure 6.2.  Three flow lines between critical points.
\end{center}

Choose a small neighborhood, $U$, of three flow lines between the two points.  Identify a neighborhood of $K_1 \cup K_2$ in $Y$, $N$, with $U$ such that each knot gets mapped to the union of two of the three flows.  We will now use $h$ to refer to the induced Morse function on $N$, with index 0 and index 3 critical points, $p$ and $q$.  Extend $h$ to a Morse function $f$ on all of $Y$ such that it is still self-indexing.  If there were no other index 0 or index 3 critical points, then we could construct the desired Heegaard diagram simply by choosing the three basepoints to be where the three flow lines pass through the Heegaard surface, $f^{-1}(3/2)$.  The idea is to cancel any critical points of index 0 or 3 outside of $N$, without affecting $f|_N$.

If such critical points exist, we rescale the Morse function in a neighborhood of $p$ and $q$ so as to not affect the critical points, but make $h(p) = -\epsilon$ and $h(q) = 3+\epsilon$ (and thus the same for $f$).   Now, remove the balls $\{f>3+\epsilon/2\}$ and $\{f<-\epsilon/2\}$ around the index 0 and index 3 critical points from $N$, to obtain a cobordism $W: S^2 \rightarrow S^2$.  In the terminology of \cite{hcob}, this is a self-indexing Morse function on the triad $(W,S^2,S^2)$.  Since each manifold in the triad is connected, we know that for each index 0 critical point, there is a corresponding index 1 with a single flow line traveling to the index 0.  This pair can be canceled such that the Morse function will not be changed outside of a neighborhood of the flow line between.  We want to see that by perhaps choosing a smaller neighborhood, $N'$, of the knots inside of $N$, this flow line does not hit $N'$.  This must be the case because if no such neighborhood existed, by compactness, this flow line would have to intersect $K_1$ or $K_2$.  But these are flows of $f$ themselves, so the two lines cannot intersect.

Hence, we can alter $f$ to remove the index 0/1 pair without affecting $f|_{N'}$.  By repeating this argument and an analogous one for index 2/3 pairs, we can remove all of the critical points of index 0 and 3 in $W$ in this fashion.  This says, after rescaling the function on the neighborhoods of $p$ and $q$ back to their original values, the new Morse function is self-indexing on $Y$ with exactly one index 0 and one index 3 critical point, and furthermore, still agrees with $h$ when restricted to a small enough neighborhood of the knots.  This is exactly what what we want to give the desired Heegaard diagram.  Now, perform an isotopy avoiding $z$, $u$, and $w$ to make the diagram admissible as necessary.
\end{proof}

\begin{remark}
This result can also be proven directly by constructing the Heegaard diagrams for each of $K_1$ and $K_2$, similar to the one associated to a knot projection in \cite{hfalternating}, but for an arbitrary $Y$.  This gives a four-pointed diagram such that proceeding to remove the correct point gives a diagram for $(K,K_1,K_2)$.
\end{remark}

Consider the link in the Kirby diagram for $M_n$, Figure 3.1.  Let the knot $Z_n$ be the one on the left which varies with respect to $n$ inside of the three-manifold obtained by 0-surgery on the remaining two components.  This is the knot that we will apply the surgery formula to.

\begin{proposition} \label{helmetexists}
For each $n$, there exists a Heegaard diagram for $(Z_1,Z_{n-1},Z_n)$ in $S^2 \times S^1 \# S^2 \times S^1$.
\end{proposition}
\begin{proof}
Let us first consider the diagram below.

\qquad\quad\includegraphics[scale=.50]{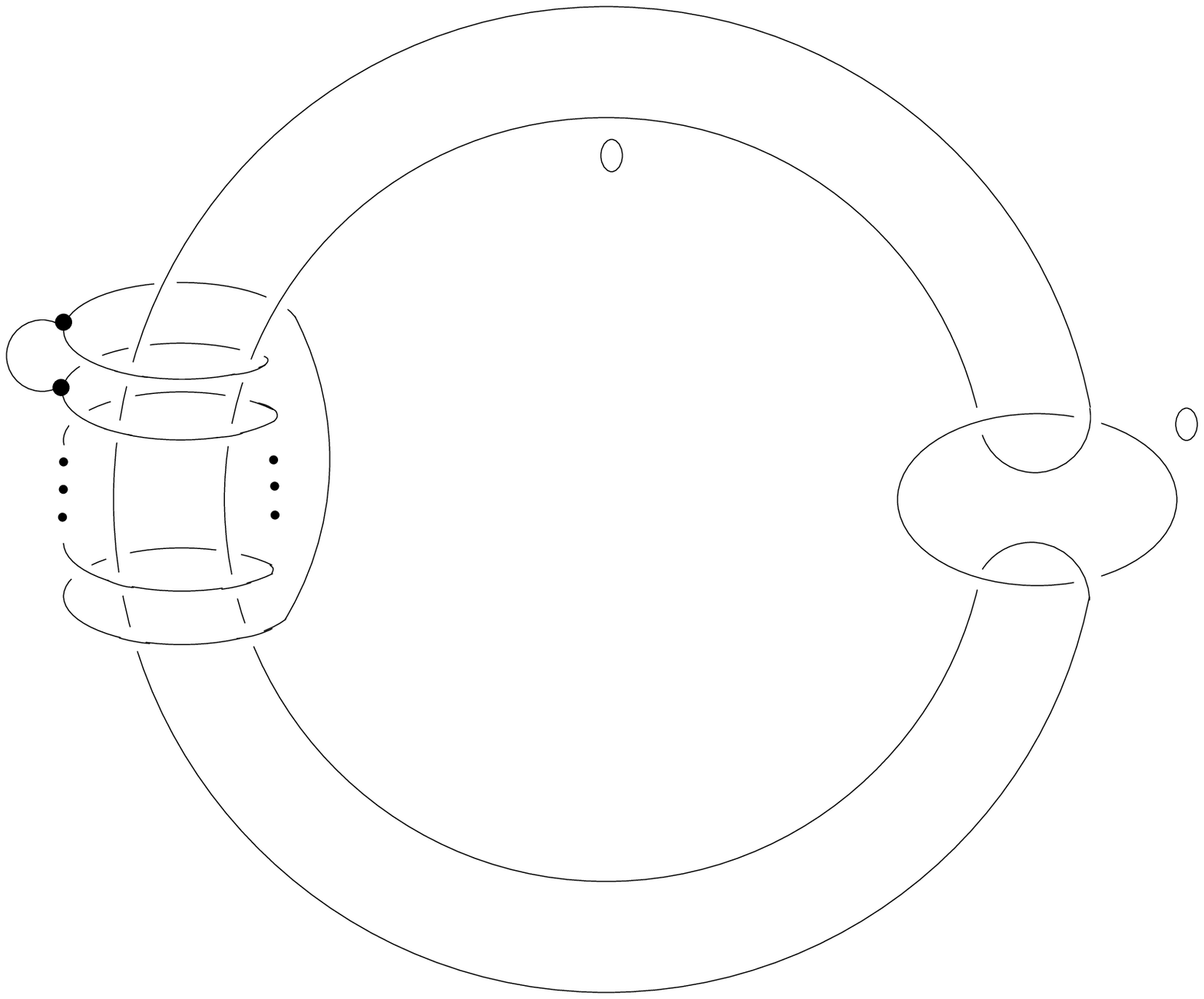}
\begin{center}
Figure 6.3.
\end{center}

Here we have attached an arc to $Z_n$ at two points (the large black
dots). This creates two additional knots as follows. Note that one
can travel two different paths from the bottom to the top attachment points; we may either wind
in an upward spiral once around the two vertical strands or follow the path that begins by
winding downward $n-1$ times. Beginning at the top attachment point,
following the attaching arc to the bottom point, and finally
traversing one of the two winding paths back to the top point gives
either $Z_1$ or $Z_{n-1}$.  Figure 6.3 now illustrates how the three knots are in position to apply
Lemma~\ref{morsehelmet}.
\end{proof}

When applying the surgery formula for $\mathbb{T}^3$, it was critical to use the map $\Theta^K$ to make all of the projections diagonal matrices.  The following lemma will allow us to do this in general.

\begin{lemma}
Consider a Heegaard diagram for $(K,K_1,K_2)$.  Then $\Theta^K = U^k \cdot \Theta^{K_2} \circ \Theta^{K_1}$ for some $k \in \mathbb{Z}$.
\end{lemma}
\begin{proof}
Consider the Alexander gradings for the three knots in the diagram.
\begin{align*}
A_K(x) - A_K(y) &= n_z(\phi) - n_w(\phi)\\ &= n_z(\phi) - n_u(\phi) + n_u(\phi) - n_w(\phi) \\ &=
A_{K_1}(x) - A_{K_1}(y) + A_{K_2}(x) - A_{K_2}(y)
\end{align*}
for each $\phi \in \pi_2(x,y)$.
\end{proof}

The following diagram will provide a useful visual reference for the upcoming proposition.

\[
\xymatrix{
\qquad && \mathcal{K}_{z,w} \ar@/_2pc/[lldd]_{\proj^{-K}} \ar@/^2pc/[ddrr]^{\proj^K}\\
&\mathcal{K}_{z,u} \ar[ld]^{\proj^{-K_1}} \ar[dr]^{\proj^{K_1}} &&\mathcal{K}_{u,w} \ar[ld]_{\proj^{-K_2}} \ar[dr]_{\proj^{K_2}} \\
\mathcal{K}_z \ar[rr]_{\destab^{-K_1}} \ar@/_2pc/[rrrr]_{\destab^{-K}}&&\mathcal{K}_u \ar[rr]_{\destab^{-K_2}} &&\mathcal{K}_w \\
}
\]

\begin{center}
Figure 6.4.
\end{center}

Fix a Heegaard diagram as given by Proposition~\ref{helmetexists}.  We now will choose the proper bases as in the $\mathbb{T}^3$ example.

\begin{proposition}\label{choosingbases}
Following Section~\ref{examplesection} for $\mathbb{T}^3$, choose bases for $\mathcal{K}_{z,u}$, $\mathcal{K}_z$, and $\mathcal{K}_u$, such that the projections and $\Theta^{Z_1}$ are given by the identity and the map $\destab^{-Z_1}$ is a matrix in $X$.  Now, choose bases for $\mathcal{K}_w$ and $\mathcal{K}_{u,w}$ such that the projections and $\Theta^{Z_{n-1}}$ are the identity.  There exists a basis for $\mathcal{K}_{z,w}$ such that $\proj^{-Z_n}$ is given by the identity, while $\proj^{Z_n}$ and $\Theta^{Z_n}$ are given by multiplication by $U^k$.
\end{proposition}

\begin{proof}
Clearly we can fix a basis for $\mathcal{K}_{z,w}$ such that $\proj^{-Z_n}$ is the identity.  Now, we combine the fact that $\proj^{Z_n} = \Theta^{Z_n} \circ \proj^{-Z_n}$ with $\Theta^{Z_n} = U^k \cdot \Theta^{Z_{n-1}} \circ \Theta^{Z_1} = U^k \cdot I$, to get the required result.
\end{proof}

\begin{remark}
By the same arguments as before, we must have that $k=0$ in the above proposition, or else $\Phi^K + \Phi^{-K}$ will have rank at least 2, contradicting the bounds coming from the spectral sequence.
\end{remark}

\begin{remark}
These constructions could be generalized to any number of basepoints (and the corresponding larger number of induced knots), but we only need three basepoints for our purposes.
\end{remark}

Although $\destab^{-Z_{n-1}}$ is not necessarily represented by an element of $X$ in this diagram, we do know that it does not contain any $U$'s in its matrix representation, since it preserves absolute grading and by construction, the basis elements of $\mathcal{K}_u$ and $\mathcal{K}_w$ in fact have the same grading.

\begin{remark}
While the individual matrix representations may seem to depend on the choice of Heegaard diagram, if $\destab^{-K} = I$, this is independent of the diagram as long as the bases are chosen such that $\proj^K = \Theta^K = I$.  A similar statement based on the work of Section~\ref{examplesection} can be made about $\destab^{-K}$ being in $X$ regardless of diagram.
\end{remark}

We are now ready for the calculation of the maps $\destab^{-Z_n}$ for all $n$.

\begin{theorem}
With this choice of bases as given by Proposition~\ref{choosingbases} for the triple $(Z_1,Z_{2n},Z_{2n+1})$, we have that $\destab^{-Z_{2n}}$ is the identity and $\destab^{-Z_{2n+1}}$ is a matrix in $X$ for all $n \geq 0$.
\end{theorem}

\begin{proof}
For $n=0$, we know that the map $\destab^{-Z_0}$ must be the identity in order to have  rk$HF^\infty(\#_{i=1}^3 S^2 \times S^1)=8$.  Similarly, from our computation for $\mathbb{T}^3$, we have seen that $\destab^{-Z_1}$ is in $X$.  Thus, the base case is established.
For the induction step, note that as soon as $\destab^{-Z_{2n}}$ is the identity, we can compose with $\destab^{-Z_1}$ to get that $\destab^{-Z_{2n+1}}$ is of type $X$.  Thus, we only need to find $\destab^{-Z_{2n}}$.  By hypothesis, $\destab^{-Z_{2n-1}} \in X$.  The key observation occurs when $\destab^{-Z_1}$ and $\destab^{-Z_{2n-1}}$ are represented by two different elements of $X$, when considering bases chosen for $(Z_1,Z_{2n-1},Z_{2n})$.  If so, then the product of the matrices, which is the representative for $\Phi^{-Z_{2n}}$, has the property that its sum with the identity, $\Phi^K$, has rank at least 2.  However, this is impossible from the spectral sequence.  Therefore, both $\destab^{-Z_{2n-1}}$ and $\destab^{-Z_1}$ are represented by the same matrix.  But, every element of $X$ squares to the identity.  $\destab^{-Z_{2n}}$ must then be the identity.
\end{proof}

We can now conclude that $HF^\infty(M_{2n},\mathfrak{s}_0)$ has rank 8 and $HF^\infty(M_{2n+1},\mathfrak{s}_0)$ has rank 6.  But, this shows exactly that $d_3$ must satisfy $x_1 \wedge x_2 \wedge x_3 \mapsto \langle x_1 \smile x_2 \smile x_3, [Y] \rangle \text{ (mod } 2)$, proving Theorem~\ref{b1=3}.

\section{Calculations for $b_1=4$}

Recall from Proposition~\ref{foursplits} that for $b_1 = 4$, $Y$ has integral triple cup product form isomorphic to that of $M_n \# S^2 \times S^1$ for some $n$ .  We then choose a basis for $H^1(Y;\mathbb{Z})$, $\{x_1, x_2, x_3, x_4\}$, with the property that $\langle x_1 \smile x_2 \smile x_3,[Y] \rangle = n$ and each $x_i$ has cup product 0 with $x_4$.

\begin{theorem}
Let $\mathfrak{s}$ be torsion.  If $n$ is even, $HF^\infty(Y,\mathfrak{s})$ has rank 16.  For $n$ odd, $HF^\infty(Y,\mathfrak{s})$ has rank 12.
\end{theorem}
\begin{proof}
As before, we simply need to calculate $HF^\infty$ for $M_n \# S^2 \times S^1$.  By the connect sum formula, $HF^\infty(Y,\mathfrak{s}) \cong HF^\infty(M_n,\mathfrak{s}_0) \otimes \mathbb{F}^2[U,U^{-1}]$.  Therefore, applying the results of the previous section gives the result.
\end{proof}

\begin{remark} Since both $M_n$ and $S^2 \times S^1$ have $HF^\infty_* \cong HF^\infty_{*+1}$, it is easy to see that this now also holds for any $Y$ with $b_1 = 4$.  These facts about grading can also be derived directly from the integer surgery formula.
\end{remark}

\begin{proof}[Proof of Corollary~\ref{b1=4}]
To see that the homology agrees with the differential coming from the conjecture, we just need to study the differential $d_3$.  If $n$ is even, then we have the result, since both homologies are rank 16, as $d_3 \equiv 0$.  Now, consider the case where $n$ is odd.
On $\Lambda^4$, $d_3: x_1 \wedge x_2 \wedge x_3 \wedge x_4 \mapsto x_4$.  This gives 3 copies of $\mathbb{F}[U,U^{-1}]$.  Now, on $\Lambda^3$, $d_3$ maps $x_1 \wedge x_2 \wedge x_3$ to $1$, and everything else to 0.  Therefore, we get 3 more copies of $\mathbb{F}[U,U^{-1}]$.  Finally, the last 6 copies of $\mathbb{F}[U,U^{-1}]$ come from the differential being 0 on $\Lambda^2$.
\end{proof}

\section{Proof of the Existence of Homologically Split Surgery Presentations}

This proof has been reproduced with the permission of Ciprian Manolescu.  \\

We start with a discussion of some results from algebra.

A {\em lattice} is a free $\Z$-module of finite rank, together with a nondegenerate symmetric bilinear form taking values in $\Z.$ A lattice $S$ is called {\em odd} if there exists $t \in S$ with $t \cdot t \in \Z$  being odd. By $S_1 \oplus S_2$ we denote the orthogonal direct sum of two lattices.

The bilinear form of a lattice $S$ determines an embedding $S \hookrightarrow S^* = \Hom(S, \Z).$ The factor group $A_S = S^*/S$ is a finite Abelian group. It comes naturally equipped with a bilinear form
$$ b_S: A_S \times A_s \to \Q/\Z, \ \ b_S(t_1 + S, t_2 + S) = t_1 \cdot t_2 + \Z,$$
called the {\em discriminant-bilinear form} of $S.$

The following results are taken from the literature; see \cite{Kneser}, \cite{Durfee}, \cite{Serre}, \cite{Nikulin}:

\begin {theorem}[Kneser-Puppe, Durfee]
\label {thm:Kneser}
Two lattices $S_1$ and $S_2$ have isomorphic discriminant-bilinear forms if and only if there exist unimodular lattices $L_1, L_2$ such that $S_1 \oplus L_1 \cong S_2 \oplus L_2.$
\end {theorem}

\begin {theorem}[Milnor]
\label {thm:Milnor}
Let $S$ be an indefinite, unimodular, odd lattice. Then $S \cong m\langle 1 \rangle \oplus n \langle -1 \rangle$ for some $m, n \geq 1.$
\end {theorem}

We say that two lattices $S_1, S_2$ are {\em stably equivalent} if there exist nonnegative integers $m_1, n_1, m_2, n_2$ such that the stabilized lattices
$$ S_1' =  S_1 \oplus m_1\langle 1 \rangle \oplus n_1\langle -1 \rangle,$$
$$ S_2' = S_1 \oplus m_2\langle 1 \rangle \oplus n_2\langle -1 \rangle$$
are isomorphic.

Note that for any lattice $S$, the direct sum $S \oplus \langle 1 \rangle \oplus \langle -1\rangle$ is indefinite and odd. Therefore, an immediate consequence of Theorems~\ref{thm:Kneser} and \ref{thm:Milnor} is:
\begin {corollary}
\label {cor:Stably}
Two lattices are stably equivalent if and only if they have isomorphic discriminant-bilinear forms.
\end {corollary}

Observe that we can restate Theorem~\ref{thm:Milnor} by saying that all unimodular lattices are stably diagonalizable. This is not the case for general lattices. Indeed, Corollary~\ref{cor:Stably} shows that a lattice is stably diagonalizable if and only if its discriminant-bilinear form comes from a diagonal lattice.
Wall \cite{Wall} classified nonsingular bilinear forms on finite Abelian groups, and showed that any such form can appear as a discriminant-bilinear form of a lattice; see also \cite[Proposition 1.8.1]{Nikulin}.
The classification contains non-diagonal forms. As a consequence, for example, the lattice of rank two given by the matrix $$H_2 = \begin{pmatrix} 0 & 2 \\ 2 & 0 \end {pmatrix}$$ is not stably diagonalizable.

Neverthless, from the classification scheme (see \cite[Proposition 1.8.2 (d)]{Nikulin}) we do obtain the following result:
\begin {proposition}
\label {prop:SD}
For any lattice $S,$ there exist a diagonal lattice $L$ (not necessarily unimodular), such that $S \oplus L$ is diagonalizable.
\end {proposition}

For example,  $H_2 \oplus \langle 2 \rangle$ is isomorphic to $\langle 2 \rangle \oplus \langle 2 \rangle \oplus \langle -2 \rangle.$

\begin {remark}
\label {rem:Sing}
Any degenerate symmetric bilinear form over $\Z$ can be expressed as a direct sum of some zeros and a non-degenerate one. Hence, the result of Proposition~\ref{prop:SD} applies to all symmetric bilinear forms (not necessarily non-degenerate).
\end {remark}

We now return to topology. Let $Y$ be a $3$-manifold. We represent it by surgery on $S^3$ along a framed link, with linking matrix $S.$ Handleslides and stabilizations correspond to elementary operations (integral changes of basis, and direct sums with $\langle \pm 1 \rangle$) on the bilinear form of $S.$  Hence, Proposition~\ref{prop:SD} and Remark~\ref{rem:Sing} complete the proof.

\end{document}